\documentclass[notitlepage,11pt]{article}
\usepackage{amssymb,mathrsfs,amsmath}
\catcode`\@=11
\@addtoreset{equation}{section}

\catcode`\@=12

\newtheorem{Theorem}{Theorem}[section]
\newtheorem{Lemma}[Theorem]{Lemma}
\newtheorem{Proposition}[Theorem]{Proposition}

\newtheorem{Remark}[Theorem]{Remark}

\newcommand{\Rn}{{\mathbb R}^{N}}
\newcommand{\Sf}{{\mathbb S}}

\newcommand{\ba} {\beta}

\newcommand{\De} {\Delta}
\newcommand{\la} {\lambda}

\newcommand{\na} {\nabla}

\newcommand{\var} {\varepsilon}

\newcommand{\R}{\mathbb{R}}
\def\proof{\noindent{\textbf{Proof. }}}
\def\QED{\hfill {$\square$}\goodbreak \medskip}

\newcommand{\un}{u_{n}}
\newcommand{\deb}{\rightharpoonup}
\newcommand{\starstar}{2^{*\!*} }

\begin{document}

\title{Entire solutions for a class of variational problems involving the biharmonic operator
and Rellich potentials}

\author{Mousomi Bhakta
\footnote
{TIFR Centre for Applicable Mathematics, Post Bag No. 6503
Sharadanagar, Chikkabommasandra, Bangalore 560065. 
Email: {mousomi@math.tifrbng.res.in}} 
\and 
Roberta Musina
\footnote
{Dipartimento di Matematica ed Informatica, Universit\`a di Udine, via delle Scienze, 206 -- 33100 Udine, Italy.
Email: {musina@uniud.it}.
Partially supported by Miur-PRIN 2009WRJ3W7-001 
{\em Fenomeni di concentrazione e problemi di analisi geometrica}
}
}

\date{}

\maketitle

\begin{abstract}
\noindent
\footnotesize We study existence, multiplicity and qualitative properties of entire
solutions for a noncompact problem related to second-order interpolation inequalities with weights. More precisely, we deal with the following family of equations
$$
{\Delta^2u=\lambda|x|^{-4}u+|x|^{-\beta}|u|^{q-2}u\quad\textrm{in $\R^N$,}}
$$
where $N\!\ge 5$, 
$q>2$, $\beta\!=\!N\!-\!q(N-4)/2$ and $\lambda\!\in\!\R$ is smaller than the Rellich constant.

%
\bigskip

\noindent
\textbf{Keywords:} {Caffarelli-Kohn-Nirenberg type inequalities, weighted biharmonic operator, dilation invariance, breaking positivity, breaking symmetry.}
\medskip

\noindent
\textit{2010 Mathematics Subject Classification:} {26D10, 47F05.}
\end{abstract}

\section{Introduction}

In this paper we study weak solutions to problem
\begin{equation}
\label{eq:problem}
\begin{cases}
\Delta^2u=\lambda|x|^{-4}u+|x|^{-\beta}|u|^{q-2}u&\textrm{in $\R^N$}\\
u\in \mathcal D^{2,2}(\R^N)~,\quad u\neq 0~.
\end{cases}
\end{equation}
Our main assumptions are the following: $N\ge 5$ and
\begin{equation}
\label{eq:ass_q_lambda}
q>2~,\quad \lambda<\gamma_N^2~,\quad \beta=N-q~\!\frac{N-4}{2}~\!,
\end{equation}
where 
\begin{equation}
\label{eq:gamma}
\gamma_N:=\frac{N(N-4)}{4}~\!.
\end{equation}
It is well known that $\gamma_N^2$ is the
best constant in the Rellich inequality
\begin{equation}\label{Rellich}
\int_{\Rn}|\De u|^2dx\geq \displaystyle\gamma_N^2\int_{\Rn}|x|^{-4}|u|^2 dx\quad\textrm{for any}\   \  u\in \mathcal D^{2,2}(\R^N),
\end{equation}
see \cite{Rel54}, \cite{Rel69}. 

Some existence and multiplicity results for (\ref{eq:problem})
in case $\lambda=0$ can be found in \cite{CM2}. If $q=\starstar$, where
$$
\starstar:=\frac{2N}{N-4}
$$
is the critical Sobolev exponent, then (\ref{eq:problem}) becomes
\begin{equation}
\label{eq:problem_critical}
\begin{cases}
\Delta^2u=\lambda|x|^{-4}u+|u|^{\starstar-2}u&\textrm{in $\R^N$}\\
u\in \mathcal D^{2,2}(\R^N)~,\quad u\neq 0.
\end{cases}
\end{equation}
In particular, if in addition $\lambda=0$ then ground state solutions to (\ref{eq:problem_critical}) 
are extremals for the Sobolev constant $S^{*\!*}$
relative to the embedding $ \mathcal D^{2,2}(\R^N)
\hookrightarrow L^{\starstar}(\R^N)$, see for instance \cite{Sw}.

The second-order version of (\ref{eq:problem_critical}), namely,
\begin{equation}
\label{eq:problem_critical_one}
\begin{cases}
-\Delta u=\lambda|x|^{-2}u+|u|^{2^*-2}u&\textrm{in $\R^N$}\\
u\in \mathcal D^{1,2}(\R^N)~,\quad u\neq 0
\end{cases}
\end{equation}
was studied by Terracini \cite{Te}, under the assumption $\lambda<(N-2)^2/4$,
the Hardy constant. We quote also \cite{CatWan}, \cite{FelSch} and references there-in
for a larger class of problems related to the Caffarelli-Kohn-Nirenberg inequalities
\cite{CKN}. 

In the first part of the paper we deal with
{\em radially symmetric} solutions to (\ref{eq:problem}). By using the Rellich inequality
(\ref{Rellich}) and the results in \cite{CM2} one can check that for any
$q,\lambda$ and $\beta$ satisfying (\ref{eq:ass_q_lambda}), the infimum
$$
S_q^{\rm rad}(\la):  =  \inf_{\scriptstyle u\in \mathcal D^{2,2}(\Rn)
\atop\scriptstyle u=u(|x|)~,~u\ne 0}\frac{\displaystyle
\int_{\Rn}|\De u|^{2}dx-\lambda\int_{\Rn}|x|^{-4}|u|^{2}dx}
{\displaystyle\left(\int_{\Rn}|x|^{-\beta}|u|^{q}dx\right)^{2/q}}
$$
is positive. The following existence result holds.

\begin{Theorem}
\label{T:existence_radial}
Let $N\ge 5$, $q>2$, $\lambda<\gamma_N^2$, and put
$\beta=N-q(N-4)/2$. Then problem (\ref{eq:problem}) has at least
one radially symmetric solution $u$ that achieves $S_q^{\rm rad}(\la)$.
\end{Theorem}

Notice that problem (\ref{eq:problem}) is invariant with respect to the
weighted dilation
\begin{equation}
\label{eq:transform}
u(x)\mapsto t^{\frac{N-4}{2}}~\!u(tx)~,\quad t>0.
\end{equation}
From now on we will identify solutions that coincide up to a weighted
dilation and change of sign. This agreement is used in the next uniqueness and positivity
result.

\begin{Theorem}
\label{T:unique}
Let $N\ge 5$, $q>2$, and put $\beta=N-q(N-4)/2$. If
$$
-(N-2)^2\le\lambda<\gamma_N^2~\!,
$$
then there exists a unique radially symmetric solution $u\in \mathcal D^{2,2}(\R^N)$ to (\ref{eq:problem}).
Moreover, $u$ is positive on $\R^N$.
\end{Theorem}

Theorems \ref{T:existence_radial} and  \ref{T:unique} will be proved in Section \ref{S:radial}. Some preliminary results 
of independent interest will be proved in Section \ref{S:ODE}. 

Notice that the positivity of $u$ does not follow
from its characterization as an extremal for $S_q^{\rm rad}(\la)$. 
Indeed, in problems involving the biharmonic operator a {\em breaking positivity}
phenomenon might appear and ground state solutions might be
forced to change sign, see for instance \cite{CM2}. Actually one could  wonder whether extremals
for $S_q^{\rm rad}(\la)$ change sign if $\lambda<<0$.
A similar phenomenon would be completely
 new with respect to (\ref{eq:problem_critical_one}) and to
similar second-order problems.

\medskip

In case $q\le \starstar$ one can use again variational methods to find
solutions to (\ref{eq:problem}) that are not necessarily radially symmetric.
We recall that
\begin{equation}\label{Sobolev}
\int_{\Rn}|\De u|^2 dx\geq S^{*\!*}\displaystyle\left(\int_{\Rn}|u|^{\starstar }dx\right)^{2/\starstar}\quad\textrm{for any}\   \  u\in \mathcal D^{2,2}(\Rn),
\end{equation}
where $S^{*\!*}$ is the Sobolev constant. If $q$ and $\beta$ are as in (\ref{eq:ass_q_lambda}), then by
interpolating (\ref{Sobolev}) and  (\ref{Rellich}) via H\"older inequality, one
plainly gets that there exists a constant $C=C(N,q)>0$ such that
\begin{equation}\label{CKN}
\int_{\Rn}|\De u|^{2}dx\ge C\left(\int_{\Rn}
|x|^{-\beta}|u|^{q}dx\right)^{2/q}\quad\textrm{for any $u\in \mathcal D^{2,2}(\Rn)$}.
\end{equation}
Notice that  (\ref{CKN}) is the second-order version of the celebrated 
Caffarelli-Kohn-Nirenberg inequalities \cite{CKN}; we quote also \cite{CM2} for a large class
of dilation-invariant inequalities on cones. It is clear that the infimum
$$
S_q(\la):  =  
\inf_{\scriptstyle u\in \mathcal D^{2,2}(\Rn)
\atop\scriptstyle u\ne 0}
\frac{\displaystyle
\int_{\Rn}|\De u|^{2}dx-\la\int_{\Rn}|x|^{-4}|u|^{2}dx}
{\displaystyle\left(\int_{\Rn}|x|^{-\beta}|u|^{q}dx\right)^{2/q}}
$$
is positive, provided that $\lambda<\gamma_N^2$.
In addition, extremals for $S_q(\la)$ give rise to solutions to (\ref{eq:problem}).
In Section \ref{S:non_radial} we
prove the next existence result.

\begin{Theorem}\label{exnonex}
Let $N\ge 5$, $q\in(2,\starstar]$, $\lambda<\gamma_N^2$, and put
$\beta=N-q(N-4)/2$. 
\begin{enumerate}
\item [(i)]
The infimum $S_q(\la)$ is achieved for any $q\in(2,\starstar)$.
\item[(ii)]
The infimum $S_{\starstar}(\la)$ is achieved if and only if
$\la\ge 0$.
\end{enumerate}
\end{Theorem}

Clearly enough, it results that $S_q(\la)\le S_q^{\rm rad}(\la)$. 
In Section \ref{S:symmetry} we wonder whether the solutions in Theorem
\ref{exnonex} are radially symmetric or breaking symmetry occurs. First, by using rearrangement techniques
we prove that $S_q(\la)=S_q^{\rm rad}(\la)$ provided that $\lambda\ge 0$.
In contrast, we show that if $\lambda<<0$ then 
$S_q(\lambda)<S_q^{\rm rad}(\la)$. Therefore, if in addition $q\in(2,\starstar)$
then problem (\ref{eq:problem}) has at least two distinct solutions.

\medskip

\small

\noindent
{\bf Notation}

We put $\R_+=(0,\infty)$ and we denote by $B_R$ the open ball in $\R^N$ of radius $R>0$ centered at $0$.

Let $\omega$ be a non-negative measurable function on 
a domain $\Omega$ in $\R^n$. For any $p\in[2,\infty)$ the weighted
Lebesgue space $L^p(\Omega;\omega(x)~\! dx)$ is the space of measurable maps
$u$ in $\Omega$ with finite norm $\left(\int_{\Omega}|u|^p
\omega(x)~\!dx\right)^{1/p}$.  If $\omega\equiv 1$ we simply write
$L^p(\Omega)$.

Let $N\ge 5$. The space $\mathcal D^{2,2}(\R^N)$ is the completion of $C^\infty_c(\R^N)$
with respect to the Hilbertian norm
$$
\|u\|^2=\int_{\R^N}|\Delta u|^2~dx.
$$
Then $\mathcal D^{2,2}(\R^N)\hookrightarrow L^2(\R^N;|x|^{-4}dx)$ by the Rellich
inequality (\ref{Rellich}), and $\mathcal D^{2,2}(\R^N)\hookrightarrow L^{\starstar}(\R^N)$
by Sobolev embedding theorem, where $\starstar=2N/(N-4)$. 
We denote by $S^{*\!*}$ the (second order) Sobolev constant, that is,
$$
S^{*\!*}=\inf_{\scriptstyle u\in \mathcal D^{2,2}(\Rn)
\atop\scriptstyle u\ne 0}
\frac{\displaystyle
\int_{\Rn}|\De u|^{2}~\!dx}
{\displaystyle\left(\int_{\Rn}|u|^{\starstar}~\!dx\right)^{2/\starstar}}~\!.
$$
It is well known that $S^{*\!*}$ is achieved by a (unique,
up to a multiplicative constant, and up to dilations and translations in
$\R^N$) positive and radially symmetric function.

If $g\in L^1(\R)$ we will often use the notation $\int g~\!ds$ instead of
$\int_{\!\!-\!\infty}^\infty g(s)~\!ds$.

\normalsize

\section{Homoclinic solutions to a fourth order ODE}
\label{S:ODE}

To prove Theorems \ref{T:existence_radial} and  \ref{T:unique} we
need some preliminary results that hold for ordinary
differential equations involving differential operators of the form
$$
Lw= w''''-2A~\!w''+B^2 w.
$$
We start with a lemma about a linear equation.

\begin{Lemma}
\label{L:unique}
Let $a\in L^\infty$, $a\ge 0$ and let $\eta\in H^2(\R_+)$ be a solution to
\begin{equation}
\label{eq:eta}
\begin{cases}
\eta''''-2A~\! \eta{''}+B^2 \eta=a(s)\eta&\textrm{on $\R_+$}\\
\eta(0)=\eta'(0)=\eta''(0)=0~\!.
\end{cases}
\end{equation}
If $A\ge B\ge 0$ then $\eta\equiv 0$ on $\R_+$.
\end{Lemma}

\proof
If $\eta'''(0)=0$ then we are done, since we immediately infer $\eta\equiv 0$.
We argue by contradiction. Replacing $\eta$ with $-\eta$ if needed, we can assume that
$\eta'''(0)>0$. Since $A\ge B$, then the equation
$c^2-2A~\!c+B^2=0$
has two real nonnegative roots $c_+, c_-$. Next we put
\begin{equation}
\label{eq:defpsi}
\psi=-\eta''+c_+\eta.
\end{equation}
Notice that in particular $\psi\in H^1(\R_+)$ and it solves
\begin{equation}
\label{eq:psi}
\begin{cases}
-\psi''+c_-\psi=a(s)\eta\quad\textrm{on $\R_+$}\\
\psi(0)=0\\
\psi'(0)<0.
\end{cases}
\end{equation}
Since 
$6\eta(s)=\eta'''(0) s^3+o(s^3)$ and $\psi(s)=-\eta'''(0)s+o(s)$  as $s\to 0^+$,
and since $\eta'''(0)>0$, then 
$\eta>0$ and $\psi<0$ in a right neighborhood of $0$. In addition, $\eta$ is smooth and vanishes at infinity. Thus 
there exists
$s_{\!\infty}\in(0,\infty]$ such that $\eta>0$ in 
$(0,s_{\!\infty})$ and $\eta\in H^1_0(0,s_{\!\infty})$. We claim that
$\psi<0$ on $(0,s_{\!\infty})$. 
If not, there exists
$s_1\in(0,s_{\!\infty})$ such that $\psi<0$ on $(0,s_1)$ and $\psi\in H^1_0(0,s_1)$. But then
from (\ref{eq:psi}) we readily get
$$
\int_0^{s_1}|\psi'|^2+c_-\int_0^{s_1}|\psi|^2=\int_0^{s_1}a(s)\eta\psi\le0~,
$$
which is impossible.
Finally, using $\eta\in H^1_0(0,s_{\!\infty})$ as test function in 
(\ref{eq:defpsi}) and recalling that 
$\eta>0$ and $\psi<0$ on $(0,s_{\!\infty})$,  we get
$$
\int_0^{s_\infty}|\eta'|^2+c_+\int_0^{s_\infty}|\eta|^2=\int_0^{s_\infty}\eta\psi\le0~,
$$
a contradiction. The Lemma is completely proved.
\QED

Next we state the main result of this section.

\begin{Theorem}
\label{T:unique_ODE}
Let $q>2$ and let $A, B$ be two given constants, such that
$A, B>0$. 
\begin{description}
\item$(i)$ The ordinary differential equation
\begin{equation}
\label{eq:ODE}
w''''-2A~\!w''+B^2w=|w|^{q-2}w\quad\textrm{on $\R$}
\end{equation}
has at least a solution $w\in H^2(\R)$. More precisely, $w$ achieves
\begin{equation}
\label{eq:I}
I(A,B):=\inf_{\scriptstyle w\in  H^2(\R)
\atop\scriptstyle w\ne 0}
\frac{\displaystyle
{\int}(|w_{ss}|^2+2A|w_{s}|^2+B^2|w|^2)ds}{\left(\displaystyle{\int}|w|^q ds\right)^{2/q}} ~\!.
\end{equation}
\item$(ii)$
If in addition
$A\ge B$ then $w$ is the unique (up to
translations, inversion $s\mapsto -s$ and change of sign) nontrivial solution to (\ref{eq:ODE}) in $H^2(\R)$. Moreover, $w$
can be taken to be even, positive and strictly decreasing on $\R_+$.
\end{description}
\end{Theorem}

\proof
Existence follows by nowadays standard variational arguments, by showing that
every minimizing sequence is
relatively compact in $H^2(\R)$ up to translations in $\R$
(use for instance the techniques in \cite{Li84}).

Next, we assume that $A\ge B$. We first notice that $w$ has 
at least one critical point, since it is smooth and vanishes at infinity. The next claim is the main step in our argument:
\begin{equation}
\label{eq:claim}
\textit{if $w'(s_0)=0$, then $w$ is even with respect to $s_0$.}
\end{equation}
To prove the claim, we first notice that we can assume that $s_0=0$. We put
$$
\eta(s):=w(s)-w(-s)~\!.
$$
We also define a measurable function $a(s)$ on $\R_+$ as follows. 
In points $s$ where $\eta(s)\neq 0$ we set
$$
a(s):= \frac{|w(s)|^{q-2}w(s)-|w(-s)|^{q-2}w(-s)}{w(s)-w(-s)}
$$
and we put  $a(s)\equiv 0$ on $\{\eta^{-1}(0)\}$. Notice that
$$
0\le a(s)\le (q-1)\max\{|w(s)|^{q-2}~\!,~\!|w(-s)|^{q-2}\}\quad\textrm{for any $s\in\R$}
$$
by the mean value theorem, so that in particular $a\in L^{\infty}(\R_+)$. Moreover,
$\eta\in H^2(\R)$ 
solves (\ref{eq:eta}), with $a$ defined as above. Thus
$\eta\equiv 0$ on $(0,\infty)$ by Lemma \ref{L:unique}, that is, $w(s)=w(-s)$ for any $s>0$. Claim (\ref{eq:claim}) is completely proved.

\medskip

\noindent
{\em Proof of the Theorem concluded. }
Since any solution $w$ to (\ref{eq:ODE}) is smooth and vanishes at infinity, we can assume 
that $w$ attains its positive maximum at $0$. Moreover, $0$ is the unique critical point
of $w$. Indeed, assume that $w'(s_0)=0$.
Then both $s\mapsto w(s)$ and 
$s\mapsto w(s_0+s)$ are even by (\ref{eq:claim}). Thus, for any $s\in \R$
we have that
$$
w(s+2s_0)=w(s_0+(s_0+s))=w(s_0-(s_0+s))=w(-s)=w(s)~\!.
$$
Hence $s_0=0$ since no function in $H^2(\R)$ can be periodic.
In particular, $w$ is strictly
positive on $\R$ and decreasing on $\R_+$.
Uniqueness then follows by the Amick and Toland theorem in  \cite{AT}. The 
proof is complete.
\QED

\begin{Remark}
\label{R:explicit}
If $A,B$ satisfy
$$
\frac{A}{B}=\frac{q^2+4}{4q}>1
$$
then the unique solution to (\ref{eq:ODE}) is given by
$$
w(s)=C\left(\cosh \nu s\right)^{-\frac{4}{q-2}},
$$
where
$$
\nu^2=\frac{(q-2)^2}{2(q^2+4)}~\!A~,
\quad C^{q-2}=\frac{2q(q+2)(3q-2)}{(q^2+4)^2}~\!A^2~\!.
$$
\end{Remark}

\begin{Remark}
\label{R:AT}
By the results in \cite{AT}, problem (\ref{eq:ODE}) has at most one homoclinic and
positive solution for any $A, B>0$. However, if $A<B$ we do not
know if there might exists nodal solutions $w\in H^2(\R)$. Actually quite
wild phenomena could appear. Consider the equation
\begin{equation}
\label{eq:AT}
w''''-2Aw''+B^2w=|w|^2\quad\textrm{on $\R$.}
\end{equation}
It is not difficult
to show that (\ref{eq:AT}) has a ground state solution $w\in H^2(\R)$ for any
$A,B>0$. Moreover, if $A\ge B$ then $w$ is positive, and hence it solves
(\ref{eq:ODE}) for $q=3$. It suffices to introduce the root 
$c_\pm>0$ to equation $c^2-2Ac+B^2=0$ as in the proof of Lemma \ref{L:unique}
and to notice that $\psi:=-w''+c_+w$ solves $-\psi''+c_-\psi=|w|^2\ge 0$. Hence
$-w''+c_+w\ge 0$, that implies $w\ge 0$.

On the other hand, it has been proved in \cite{CT} that a { plethora}
(countably infinite set)
of sign-changing homoclinic solutions exists for $A/B<1$ close enough to $1$.
\end{Remark}

\begin{Remark}
\label{R:Hamilton}
The differential equation in (\ref{eq:ODE}) is conservative, with Hamiltonian
$$
H(w,w',w'',w''')=
-w'''w'+\frac{1}{2}|w''|^2+A|w'|^2-\frac{1}{2}B^2|w|^2+
\frac{1}{q}|w|^q~\!.
$$
Therefore, any homoclinic solution $w$ to (\ref{eq:ODE}) satisfies
\begin{equation}
\label{eq:Hamilton}
-w'''w'+\frac{1}{2}|w''|^2+A|w'|^2-\frac{1}{2}B^2|w|^2+
\frac{1}{q}|w|^q\equiv 0\quad\textit{on $\R$.}
\end{equation}
Using (\ref{eq:Hamilton}) one infers the following a-priori
bound on homoclinic solutions to (\ref{eq:ODE}):
$$
\|w\|_\infty^{q-2}\le \frac{q}{2}~\!B^2~\!.
$$
\end{Remark}

\section{Radially symmetric solutions}
\label{S:radial}

In this section we prove Theorems \ref{T:existence_radial} and
\ref{T:unique}. The main tool is the Emden-Fowler transform,
which has already been used in \cite{CM2}.

Given a radially symmetric function $u$ we define the Emden Fowler transform of $u$ by $w=Tu$, where
$$
u(x)=|x|^{-\frac{N-4}{2}}w(-\log |x|).
$$
It turns out that
$T$ maps radially symmetric functions $u\in D^{2,2}(\Rn)$ into functions $w
\in H^2(\R)$, and moreover
\begin{gather}
\label{eq:EF1}
\int_{\Rn} |x|^{-\ba}|u|^q dx= \omega_{\!N}{\int}|w|^q ds\\
\label{eq:EF2}
\int_{\Rn}|\De u|^2 dx= \omega_{\!N}{\int}(|w_{ss}|^2+2(\gamma_N+2)|w_{s}|^2+{\gamma_N}^{2}|w|^2)ds,
\end{gather}
where $\omega_{\!N}=|\mathbb S^{N-1}|$ and
$\gamma_N=N(N-4)/4$, as in (\ref{eq:gamma}).

\bigskip

\noindent
{\bf Proof of Theorem \ref{T:existence_radial}.}
From (\ref{eq:EF1}) and (\ref{eq:EF2}) it immediately follows that
$u$ achieves $S^{\rm rad}_{q}(\la)$ if and only if $w=Tu$ attains the infimum
$I(\gamma_N+2,\gamma_N^2-\lambda)$, see (\ref{eq:I}). The conclusion follows from
$(i)$ of Theorem  \ref{T:unique_ODE}.
\QED

\bigskip

\noindent
{\bf Proof of Theorem \ref{T:unique}.}
A direct computation
shows that $u\in  \mathcal D^{2,2}(\Rn)$ is a radial solution to (\ref{eq:problem})
if and only if 
$Tu=w\in H^2(\R)$ is a solution to 
$$
w''''-2(\gamma_N+2)w''+(\gamma_N^2-\lambda)w=|w|^{q-2}w\quad\textrm{on $\R$}.
$$
Since $\lambda\ge-(N-2)^2$, then
$$
(\gamma_N+2)^2-(\gamma_N^2-\lambda)=4(\gamma_N+1)+\lambda=(N-2)^2+\lambda\ge 0~\!.
$$
The conclusion readily follows from the second part of Theorem
\ref{T:unique_ODE}.
\QED

\begin{Remark}
Using Remark \ref{R:explicit}, one can explicitly compute the
radially symmetric solution to (\ref{eq:problem}) when
$$
\lambda=\lambda(q):=\gamma_N^2-\left(\frac{4q}{q^2+4}\right)^2
(\gamma_N^2+2)^2~\!.
$$
Notice that $\lambda(q)>-(N-2)^2$ for any $q\in(2,\starstar]$ and that
$\lambda(\starstar)=0$.
If $\lambda=\lambda(q)$ then the unique radially symmetric solution to (\ref{eq:problem})
is given by
$$
u(x)=\widetilde C~\!|x|^{\frac{4-N}{2}+\frac{4\nu}{q-2}}~\!
\left(1+|x|^{2\nu}\right)^{-\frac{4}{q-2}}
$$
for a computable constant $\widetilde C$ depending on $N, q$, where
$$
\nu=\frac{(q-2)^2}{2(q^2+4)}~\!(\gamma_N^2+2)~\!.
$$
\end{Remark}

The preliminary results in Section
\ref{S:ODE} can be applied to a larger class of fourth order equations in $\R^N$.
Let us consider for instance a problem taken from \cite{CM2}.

Let $N\ge 2$ be any integer, and let $\alpha$ be a given exponent, such that
$$
\alpha>4-N~,\quad \alpha\neq N~\!.
$$
By the results in \cite{CM1}, \cite{CM2}, it turns out that for any $q\in[2,\infty)$ the infimum
$$
S^{\rm rad}_{q,\alpha}:= \inf_{\scriptstyle u\in 
\mathcal N^2_{\rm rad}(\R^N;\alpha)\atop\scriptstyle u\ne 0}
\frac{\displaystyle
\int_{\R^{n}}|x|^{\alpha}|\Delta u|^{2}dx}
{\displaystyle\left(\int_{\R^{n}}|x|^{-\beta_\alpha}|u|^{q}dx\right)^{2/q}}
$$
is positive. Here 
$$
\beta_\alpha:=N-q~\!\frac{N-4+\alpha}{2}
$$
and $\mathcal N^2_{\rm rad}(\R^N;\alpha)$ is the completion of the space
of radially symmetric functions in $C^\infty_c(\R^N)$ with respect to the
Hilbertian norm
$$
\|u\|^2=\int_{\R^N}|x|^{\alpha}|\Delta u|^2~\!dx.
$$
More precisely, for $q=2$ the {\em radial Rellich constant} is given by
$$
S^{\rm rad}_{2,\alpha}=\gamma_{N\!,\alpha}^2~,\quad
\textrm{where}~~
\gamma_{N\!,\alpha}=\frac{(N-4+\alpha)(N-\alpha)}{4}~\!,
$$
and it is not achieved. In Theorem 2.8 of \cite{CM2} it is
proved that, for any $q>2$, the infimum  $S^{\rm rad}_{2,\alpha}$ is attained 
by some $u\in\mathcal N^2_{\rm rad}(\R^N;\alpha)$. Then, up to a Lagrange
multiplier, $u$ solves
\begin{equation}
\label{eq:problem_alpha}
\Delta\left(|x|^\alpha\Delta u\right)=
|x|^{-\beta_\alpha}|u|^{q-2}u\quad\textrm{on $\R^N$.}
\end{equation}

We point out the following positivity and uniqueness result.

\begin{Theorem}
Let $\alpha>4-N$, $\alpha\neq N$. Then problem (\ref{eq:problem_alpha})
has a unique nontrivial solution $u\in \mathcal N^2_{\rm rad}(\R^N;\alpha)$.
Moreover, $u> 0$ on $\R^N$.
\end{Theorem}

\proof
Existence follows from Theorem 2.8 in \cite{CM2}. Following \cite{CM1}, we 
introduce the value
$$
\overline\gamma_{N\!,\alpha}=
\left(\frac{N-2}{2}\right)^2+\left(\frac{\alpha-2}{2}\right)^2.
$$
By direct computations (see also \cite{CM1}, \cite{CM2}), we have that
$u\in \mathcal N^2_{\rm rad}(\R^N;\alpha)$ solves (\ref{eq:problem_alpha}) if and
only if $w\in H^2(\R)$ solves
$$
w''''-2\overline\gamma_{N\!,\alpha}~\!w''+
\gamma_{N\!,\alpha}^2w=|w|^{q-2}w\quad\textrm{on $\R$.}
$$
Here we have denoted by $w=T_\alpha u$ the weighted Emden-Fowler transform
of $u$:
$$
u(x)=|x|^{\frac{4-N-\alpha}{2}}~\!w(-\log|x|).
$$
Since $\overline\gamma_{N\!,\alpha}\ge\gamma_{N\!,\alpha}$, then
positivity and uniqueness follows by Theorem \ref{T:unique_ODE}.
\QED

\section{On the infimum $S_{q}(\la)$}
\label{S:non_radial}

The following {\em $\var-$compactness} result is a
basic tool in the proof of Theorem \ref{exnonex}.

\begin{Proposition}\label{prop1}
Let $\un\in\mathcal D^{2,2}(\Rn)$ such that $\un\deb 0$ in $\mathcal D^{2,2}(\Rn)$ and 
\begin{equation}\label{assump1}
\De^2\un-\la |x|^{-4}\un=|x|^{-\ba}|\un|^{q-2}\un+f_n
\end{equation}
\begin{equation}\label{assump2}
\int_{B_{R}}|x|^{-\ba}|\un|^q dx\leq \var_{0}\quad\textit{for some $\var_0, R>0$,}
\end{equation}
where $f_n\to 0$ in  
the dual space of $\mathcal D^{2,2}(\Rn)$. If $\var_{0}<S_{q}(\la)^{\frac{q}{q-2}}$,
then 
$$
|x|^{-\ba}|\un|^q\to 0\quad\textit{in $L^1_{\rm loc}(B_R)$.}
$$
\end{Proposition}

\proof
Fix any $R'\in (0, R)$ and choose a cut off function $\phi\in C^{\infty}_{c}(B_R)$
such that $\phi\equiv 1$ in $B_{R'}$. Since $\na\phi$ has compact support in $\Rn\setminus\{0\}$,
then
 $$
 \int_{\Rn}\De\un\De(\phi^2\un)~\!dx= \int_{\Rn}|\De(\phi\un)|^2~\!dx+o(1)
 $$
 by Rellich Theorem.
Applying \eqref{assump1}, H\"{o}lder inequality and (\ref{assump2}) we obtain 
 \begin{eqnarray*}
\int_{\Rn}[|\De(\phi\un)|^2-\la|x|^{-4}|\phi\un|^2]dx &=&\int_{\Rn}|x|^{-\ba}|\un|^{q-2}|\phi\un|^2dx
 \\
&\leq& \var_{0}^\frac{q-2}{q}\displaystyle\left (\int_{\Rn}|x|^{-\ba}|\phi\un|^q dx\right)^\frac{2}{q}~\!.
\end{eqnarray*}
We estimate from below the left hand side of the above inequality by 
$$
 \int_{\Rn}[|\De(\phi\un)|^2-\la|x|^{-4}|\phi\un|^2]dx \ge
 S_{q}(\lambda) \left(\int_{\Rn}|x|^{-\ba}|\phi\un|^q dx\right)^{2/q}.
 $$
  Thus
\begin{displaymath}
S_{q}(\lambda)\displaystyle\left (\int_{\Rn}|x|^{-\ba}|\phi\un|^q dx\right)^{2/q}\leq \var_{0}^\frac{q-2}{q}\displaystyle\left (\int_{\Rn}|x|^{-\ba}|\phi\un|^q dx\right)^{2/q}~\!.
\end{displaymath}
Since $\var_{0}^\frac{q-2}{q}<S_q(\lambda)$ and $\phi\equiv 1$ on $B_{R'}$,
the conclusion follows by the arbitrariness of $R'\in(0,R)$.
\QED

\subsection{Proof of Theorem \ref{exnonex} in case $q<\starstar$}

Using Ekeland's variational principle we can fix a minimizing sequence for
$S_q(\lambda)$ such that
\begin{eqnarray}\label{assump4}
\int_{\Rn}|\De\un|^2-\la\int_{\Rn}{|x|^{-4}}|\un|^2 &=& \int_{\Rn}|x|^{-\ba}|\un|^q dx+o(1)\nonumber
\\
&=& S_{q}(\la)^\frac{q}{q-2}+o(1)~\!,
\end{eqnarray}
\begin{equation}\label{assump3}
\Delta^2\un-\la |x|^{-4}\un=|x|^{-\ba}|\un|^{q-2}\un+f_n~\!,
\end{equation}
where $f_n\to 0$ in the dual of $\mathcal D^{2,2}(\Rn)$. Up to a rescaling, we can also assume that
\begin{equation}\label{a}
\int_{B_{2}}|x|^{-\ba}|\un|^q dx=\frac{1}{2}S_{q}(\la)^\frac{q}{q-2}~\!.
\end{equation}
Finally, since $\la<\gamma_N^2$, then 
$\un$ is a bounded sequence in
$\mathcal D^{2,2}(\Rn)$ by Rellich inequality. Therefore we can assume that 
there exists $u\in \mathcal D^{2,2}(\Rn)$ such that $\un\deb u$ 
weakly in
$\mathcal D^{2,2}(\Rn)$. By standard arguments, to conclude it suffices to show that $u\neq 0$.
We argue by contradiction. Suppose $\un\deb 0$. 
We
can use Proposition \eqref{prop1} to get
$$
o(1)=\int_{B_{1}}|x|^{-\ba}|\un|^q dx=\int_{B_{2}}|x|^{-\ba}|\un|^q dx-
\int_{1<|x|<2}|x|^{-\ba}|\un|^q dx~\!.
$$
Thus from  (\ref{a}) we infer
\begin{equation}
\label{eq:new}
\int_{1<|x|<2}|x|^{-\ba}|\un|^q dx=\frac{1}{2}S_{q}(\la)^\frac{q}{q-2}+o(1)
\end{equation}
which readily leads to a contradiction by Rellich theorem, as $q\in(2,\starstar)$. Proposition $(i)$ in Theorem \ref{exnonex}
is completely proved.
\QED

\subsection{Proof of Theorem \ref{exnonex} in the limiting case $q=\starstar$}
We start by pointing out a sufficient condition for existence.

\begin{Lemma}\label{prop2}
If $S_{\starstar}(\la)<S^{*\!*}$ then $S_{\starstar}$ is achieved.
\end{Lemma}

\proof
As in the proof of part (i) we select a minimizing sequence for
$S_{\starstar}(\la)$ satisfying
\begin{eqnarray*}
\int_{\Rn}|\De\un|^2~\!dx-\la\int_{\R^N}{|x|^{-4}}|\un|^2~\!dx&=& \int_{\Rn}|\un|^{\starstar} dx+o(1)\\
&=&
 S_{\starstar}(\la)^{N/4}+o(1)
 \end{eqnarray*}
 \begin{equation}\label{assump6}
\Delta^2\un-\la |x|^{-4}\un=|\un|^{\starstar-2}\un+f_n
\end{equation}
$$
\int_{B_{2}}|\un|^{\starstar} dx = \frac{1}{2} S_{\starstar}(\la)^{N/4}~\!,
$$
where $f_n\to 0$ in the dual of $\mathcal D^{2,2}(\Rn)$. In addition, we can assume that
$\un$ converges weakly in $\mathcal D^{2,2}(\Rn)$. Again, we only have to exclude that 
the weak limit vanishes. By contradiction, assume that
$\un\deb 0$ in $\mathcal D^{2,2}(\Rn)$. 
Arguing as in the proof of part $(i)$ we can conclude that (\ref{eq:new}) holds with
$q=\starstar$, that is,
\begin{equation}\label{b}
\int_{1<|x|<2}|\un|^{\starstar}~\!dx=\frac{1}{2}S_{\starstar}(\la)^{N/4}+o(1).
\end{equation}
Now we choose $\phi$ to be a cut off function in $C_{c}^{\infty}(\Rn\setminus\{0\})$ such that $\phi\equiv 1$ in $B_2\setminus B_1$. We use $\phi^2\un$ as a test function in \eqref{assump6}. Rellich 
theorem and H\"older inequality plainly lead to
$$
\int_{\Rn}[|\De(\phi\un)|^2-\la|x|^{-4}|\phi\un|^2]dx \leq S_{\starstar}(\la)\displaystyle\left(\int_{\Rn}|\phi\un|^{\starstar}dx\right)^\frac{N-4}{N}+o(1)
$$
(argue as in the proof of part $(i)$). Since $\phi$ has compact support in
$\R^N\setminus\{0\}$, then Rellich theorem again and the Sobolev inequality give
$$
\int_{\Rn}[|\De(\phi\un)|^2-\la|x|^{-4}|\phi\un|^2]=
\int_{\Rn}|\De(\phi\un)|^2+o(1) \ge
S^{*\!*}\displaystyle\left(\int_{\Rn}|\phi\un|^{\starstar}\right)^{2/\starstar}.
$$
In conclusion, we have proved that
$$
S^{*\!*}\displaystyle\left(\int_{\Rn}|\phi\un|^{\starstar}dx\right)^{2/\starstar}\leq S_{\starstar}(\la)\displaystyle
\left(\int_{\Rn} {\Rn}|\phi\un|^{\starstar}dx\right)^{2/\starstar}~\!,
$$ 
which implies $\int_{\Rn}|\phi\un|^{\starstar} =o(1)$.
Hence $\int_{B_2\setminus B_1}|\un|^{\starstar}dx=o(1)$, since
$\phi\equiv 1$ in $B_2\setminus B_1$ This is a contradiction to \eqref{b}. Thus $u\not=0$
and achieves $S_{\starstar}(\la)$.
\QED

\medskip
\noindent
{\bf Proof of Theorem \ref{exnonex} - $\bf (ii)$.}
We start by recalling that the Sobolev constant $S^{*\!*}$ is achieved
in $\mathcal D^{2,2}(\R^N)$, see for instance \cite{Sw}. Next, let
$u$ be any fixed function in $C^\infty_c(\R^N\setminus\{0\})$, and put
$v_y(x)=u(x+y)$. Since
\begin{eqnarray*}
S_{\starstar}(\la) &\leq& \lim_{|y|\to\infty}\frac{\displaystyle\int [|\De v_y|^2-\la|x|^{-4}|v_y|^2]~\!}{\displaystyle\left(\displaystyle\int_{\Rn}|v_y|^{\starstar}~\!\right)^{2/\starstar}}
= \lim_{|y|\to\infty}\frac{\displaystyle\int [|\De u|^2-\la|x-y|^{-4}|u|^2]~\!}{\displaystyle\left(\displaystyle\int_{\Rn}|u|^{\starstar}~\!\right)^{2/\starstar}}
\\
&=&
\frac{\displaystyle\int |\De u|^2~\!}{\displaystyle\left(\displaystyle\int_{\Rn}|u|^{\starstar}~\!\right)^{2/\starstar}},
\end{eqnarray*}
then $S_{\starstar}(\la)\leq S^{*\!*}$ for every $\la\in\R$. Now assume $\la<0$.
Then $S_{\starstar}(\la)\geq S^{*\!*}$, that is, $S_{\starstar}(\la)=S^{*\!*}$. 
Therefore $S_{\starstar}(\la)$ can not be achieved since $S^{*\!*}$ is achieved.

Now assume $0<\la<\gamma_N^2$ and let $U$ be an extremal for $S^{*\!*}$. Then
$$
S_{\starstar}(\la) \leq\frac{\displaystyle\int_{\Rn} [|\De U|^2-\la|x|^{-4}|U|^2]~\!}{\displaystyle\left(\displaystyle\int_{\Rn}|U|^{\starstar}~\!\right)^{2/\starstar}}
< \frac{\displaystyle\int_{\Rn}|\De U|^2~\!}{\left(\displaystyle\int_{\Rn}|U|^{\starstar}~\!\right)^{2/\starstar}}= S^{*\!*}.
$$
and hence $S_{\starstar}(\la)$ is achieved
 by Lemma \ref{prop2}. 
\QED

\section{Positivity, symmetry and breaking symmetry}
\label{S:symmetry}

In this section we 
study the symmetry and breaking symmetry of $S_{q}(\la)$ depending on the parameter $\la$. 
We
identify functions $u$ that coincide
up to a multiplicative constant and up to a transform of the kind (\ref{eq:transform}).

We start by recalling that in case $\lambda=0$ and $q=\starstar$, the Sobolev
constant $S^{*\!*}=S_\starstar(0)$ is achieved by the radially symmetric function
$$
U(x)=\left(1+|x|^2\right)^{\frac{4-N}{2}}~\!,
$$
see for instance \cite{Sw}. Moreover, $u$ achieves $S^{*\!*}$ if and only if $u$ is the composition
of $U$ with a translation in $\R^N$ (here we use the above identification).
Notice that in particular $U$ is positive. 

Since truncations $u\mapsto u^\pm$ are not allowed in dealing with fourth order differential operators, 
the positivity of extremals for $S_q(\lambda)$ does not follow by standard arguments.

In the first result we use rearrangement techniques and the uniqueness result
in Theorem \ref{T:unique} to investigate the case $\lambda\ge 0$.

\begin{Theorem}
\label{T:symmetry}
Assume $\lambda\neq 0$ or $q<\starstar$.
If $\la\geq 0$ then $S_{q}(\la)$ is achieved by a unique function $u\in\mathcal D^{2,2}(\R^N)$. Moreover, $u$ is positive, radially symmetric about the origin
and radially decreasing.
\end{Theorem}

\proof
Let $u$ to be an extremal of $S_{q}(\la)$ and 
denote by  $(-\De u)^{*}$ the Schwartz symmetrization of $-\De u$. We use Lemma \ref{L:regularity} in the appendix to introduce 
 the function $v\in\mathcal D^{2,2}(\Rn)$ satisfying 
$$-\De v=(-\De u)^{*}~\!.$$ 
In turns out that $u^{*}\leq v$ on $\R^N$, 
see for instance Remark II.13 in  \cite{Li84}. Clearly, if 
 $u=u^{*}$ then we are done. Assume by contradiction that 
 $u\neq u^*$. By the theory of symmetrisation (see Lieb and Loss \cite{LL}, Theorem 3.4),
we first obtain
$$
\int_{\Rn}|\De v|^2~\!dx=\int_{\Rn}|(-\De u)^{*}|^2~\!dx=\int_{\Rn}|\De u|^2~\!dx~.
$$
In addition, since we are assuming that $u^*\neq u$, then
$$
\int_{\Rn}|x|^{-4}|u|^2~\!dx < \int_{\Rn}|x|^{-4}(|u|^2)^*dx= \int_{\Rn}|x|^{-4}|u^*|^2dx\leq \int_{\Rn}|x|^{-4}|v|^2dx~\!.
$$
Thus we infer that
$$
\lambda\int_{\Rn}|x|^{-4}|u|^q~\!dx\le\lambda\int_{\Rn}|x|^{-4}|v|^q~\!dx~\!,
$$
and that the strict inequality holds if $\lambda>0$.
Similarly, we find
$$
\int_{\Rn}|x|^{-\beta}|u|^q~\!dx\le \int_{\Rn}|x|^{-\beta}|v|^q~\!dx~\!,
$$
and the strict inequality holds if $\beta>0$, that is, if $q<\starstar$.
In conclusion, since we are assuming that $\lambda$ and $\beta$ are not
contemporarily zero, we have that
$$
S_{q}(\la)\le
\frac{\displaystyle \int_{\Rn}[|\De v|^2-\la|x|^{-4}|v|^2]dx}{
\left(\displaystyle \int_{\Rn}|x|^{-\beta}|v|^qdx\right)^{2/q}}\\
<
\frac{\displaystyle \int_{\Rn}[|\De u|^2-\la|x|^{-4}|u|^2]dx}{
\left(\displaystyle \int_{\Rn}|x|^{-\beta}|u|^qdx\right)^{2/q}}=S_{q}(\la)~\!,
$$
a contradiction.
Therefore $u=u^{*}$, that is, $u$ is nonnegative and radially symmetric decreasing function. 

Uniqueness follows immediately by using the Emden-Fowler transform
as in Section \ref{S:radial} and Theorem \ref{T:unique}.
\QED

As soon as $\lambda\to -\infty$, a braking symmetry phenomenon appears. The next theorem
applies in case $q<\starstar$, due to the
nonexistence result pointed out in the critical case $q=\starstar$,
$\lambda<0$. We quote \cite{FelSch}, and \cite{DELT} for remarkable
breaking symmetry results for similar second-order equations.

\begin{Theorem}
If $\la<<0$ then $S_{q}(\la)<S_{q}^{\rm rad}(\la)$ and hence no extremal for $S_{q}(\la)$ is radially symmetric. 
\end{Theorem}

\proof
We already know that $S_{q}(\la)\leq S_{q}^{\rm rad}(\la)$. We will give 
an explicit condition
on $\lambda$ to have $S_{q}(\la)<S_{q}^{\rm rad}(\la)$.
Define 
$$
n(u)=\int_{\Rn}[|\De u|^2-\la|x|^{-4}|u|^2]~\!dx~,\quad d(u)= \left(\int_{\Rn}|x|^{-\ba}|u|^q~\!dx\right)^{2/q}
$$ 
and $Q(u)=n(u)/d(u)$. 
Assume that 
$u$ is a radially symmetric minimizer of $Q$ on $\mathcal D^{2,2}(\Rn)$. Our goal is to show that
$-\lambda$ can not be too large.  By homogeneity we can assume that $d(u)=1$. Thus 
$Q{'}(u)\cdot v=0$ and $ Q{''}(u)[v, v]\geq 0$ for all $v\in\mathcal D^{2,2}(\Rn)$, that is,
$n{'}(u)\cdot v=Q(u) d{'}(u)\cdot v$ and
\begin{equation}
\label {c}
n{''}(u)[v, v]\geq Q(u) d{''}(u)[v, v]
\end{equation}
for all $v\in\mathcal D^{2,2}(\Rn)$.
We compute $n{''}(u)[v, v]=2\displaystyle\int_{\Rn}\left(|\De v|^2-\la|x|^{-4}|v|^2\right)~\!$ and 
$$
d{''}(u)[v, v]=2(2-q)\displaystyle\left(\int_{\Rn}|x|^{-\ba}|u|^{q-2}uv ~\!\right)^2
+2(q-1)\int_{\Rn}|x|^{-\ba}|u|^{q-2}v^2 ~\!.
$$
Now we choose the test function $v$. Let $\varphi_1\in H^{1}(\mathbb S^{N-1})$ be an eigenfunction of Laplace-Beltrami operator on $\mathbb S^{N-1}$ corresponding to the smallest positive eigenvalue. Thus
$$-\De_{\sigma}\varphi_1=(N-1)\varphi_1~, \   \  \frac{1}{|\mathbb S^{N-1}|}\int_{\mathbb S^{N-1}}|\varphi_1|^{2}~\!d\sigma=1~, \   \  \int_{\mathbb S^{N-1}}\varphi_1~\!d\sigma=0.$$
Define the test function $v$ as $v=u\varphi_1.$ Thus
$$
d{''}(u)[v, v]=2(q-1)\int_{\Rn}|x|^{-\ba}|u|^{q}=2(q-1)
$$
and
\begin{eqnarray*}
\int_{\Rn}|\De v|^2 ~\! &=& \int_{\Rn}|\De(u\varphi_1)|^2 ~\!= 
\int_{\Rn}|\De u-(N-1)|x|^{-2}u|^2\\
&\le&\int_{\Rn}[|\De u|^2+(N-1)^2|x|^{-4}u^2]
\\&& \quad\quad +
2(N-1)\displaystyle\left(\int_{\Rn}|x|^{-4}{u^2}~\!\right)^{1/2}\displaystyle\left(\int_{\Rn}|\De u|^2 ~\!\right)^{1/2}
\end{eqnarray*}
by H\"older inequality.
Therefore from \eqref{c} and from the definition of $Q(u)=n(u)$
we obtain
\begin{eqnarray*}
\nonumber
(q-2) \int_{\Rn}|\De u|^2 ~\!&\le& \left((N-1)^2+\lambda(q-2)\right)\int_{\R^N}|x|^{-4}u^2~\!~\!
\\
&&\quad\quad +
2(N-1)\displaystyle\left(\int_{\Rn}|x|^{-4}{u^2}~\!\right)^{1/2}\displaystyle\left(\int_{\Rn}|\De u|^2 ~\!\right)^{1/2}.
\end{eqnarray*}
In particular, the quantity
$$
X:=\left(\frac{\displaystyle\int_{\Rn}|\De u|^2 ~\!}
{\displaystyle\int_{\Rn}|x|^{-4}u^2~\!}\right)^{1/2}
$$
satisfies the inequality
$$
(q-2)X^2-2(N-1)X-\left((N-1)^2+\lambda(q-2)\right)\le 0~\!,
$$
and hence it necessarily holds that
$$
\lambda\ge \min_{t\in\R}
\left(~\!t^2-\frac{2(N-1)}{q-2}t-\frac{(N-1)^2}{q-2}\right)=-\frac{(q-1)(N-1)^2}{(q-2)^2}~\!.
$$
Thus, if
$$
\lambda< -\frac{(q-1)(N-1)^2}{(q-2)^2}~\!.
$$
then no extremal for $S_{q}(\la)$ is radially symmetric and  symmetry breaking occurs.
\QED

\begin{Remark}
If $q$ is close enough to $\starstar$ then one can obtain a better estimate
on the breaking symmetry parameter $\lambda$ by arguing as follows.
Notice that 
$X> \gamma_N$ by the Rellich inequality  (\ref{Rellich}). Thus, if
$\gamma_N\ge \frac{N-1}{q-2}$ that is, if
$$
2+\frac{4(N-1)}{N(N-4)}\le q\le\starstar~\!,
$$
then the radial solution $u$ does not achieve $S_q(\lambda)$ 
unless
$$
\lambda> \min_{t\ge \gamma_N}
\left(~\!t^2-\frac{2(N-1)}{q-2}t-\frac{(N-1)^2}{q-2}\right)=
 \gamma_N^2-\frac{N-1}{q-2}~\!(N-1+2\gamma_N)~\!.
$$
Conversely, if 
$$
\lambda\le \gamma_N^2-\frac{N-1}{q-2}~\!(N-1+2\gamma_N)
$$
then breaking symmetry occurs.
\end{Remark}

\appendix

\section{\!\!\!\!\!\!ppendix} 

The following lemma has been used in the proof Theorem
\ref{T:symmetry} (see also  \cite{Li84}, Remark II.13).

\begin{Lemma}
\label{L:regularity}
Let $N\ge 5$ and $g\in L^2(\Rn)$. Then the equation
$$
-\Delta v=g\quad\textit{on $\Rn$}
$$
has a unique solution $v\in \mathcal D^{2,2}(\Rn)$.
\end{Lemma}

Due to the unboundedness of the domain, Lemma \ref{L:regularity} does not follows by using a direct  variational approach. Indeed, in general,
the solution $v$ in Lemma \ref{L:regularity} does not belong to
$\mathcal D^{1,2}(\R^N)$. For instance,
for any $t\in [2,4)$ the smooth function
$$
v(x)=\left(1+|x|^2\right)^{\frac{t-N}{4}}
$$
satisfies $\Delta v\in L^2(\R^N)$ but $\displaystyle\int_{\R^N}|\nabla v|^2~\!dx=\infty$.

\medskip

Lemma \ref{L:regularity} can be proved by
investigating the integrability properties of the Green's function
at infinity. However, we will use here an alternative approach 
that allows us to underline the unexpected role of the Rellich inequality.

\medskip

Actually, to better understand the phenomena 
involved we will consider a larger class of non-homogeneous
linear problems on cones in $\R^N$, $N\ge 2$. We emphasize the fact that
low dimensions are included. We first introduce some notation. 

Let $\Sigma$ be a domain of class $C^2$ in the unit sphere $\Sf^{N-1}$. We
introduce the {\em cone}
$$
\mathcal{C}_{\Sigma}:=\left\{~r\sigma~|~r>0~,~\sigma\in\Sigma~\right\}.
$$
We denote by $\lambda_\Sigma$ the bottom
of the Dirichlet spectrum of the Laplace-Beltrami operator on $\Sigma$. For instance,
if $\Sigma=\Sf^{N-1}$ then $\mathcal C_{\Sf^{N-1}}=\R^N\setminus\{0\}$ and
$\lambda_{\Sf^{N-1}}=0$. If $\Sigma$ is has compact closure in $\Sf^{N-1}$ then
$\lambda_\Sigma>0$.
Then the following Hardy inequality holds.

\begin{Lemma}
\label{L:Hardy}
Let $\Sigma$ be a domain of class $C^2$ in the unit sphere $\Sf^{N-1}$, with
$N\ge 2$. Then
\begin{equation}
\label{eq:Hardy}
\int_{\mathcal C_\Sigma}|x|^{-2}|\nabla u|^2~dx\ge
\left[\left(\frac{N-4}{2}\right)^2+\lambda_\Sigma\right]
\int_{\mathcal C_\Sigma}|x|^{-4}|u|^2~dx
\end{equation}
for any $u\in C^\infty_c(\mathcal C_\Sigma)$.
\end{Lemma}

\proof
Fix any $u\in C^\infty_c(\mathcal C_\Sigma)$. Use the divergence theorem
and Proposition 1.1 in \cite{FM1} to estimate
\begin{eqnarray*}
\int_{\mathcal C_\Sigma}|x|^{-2}|\nabla u|^2~dx&=&
\int_{\mathcal C_\Sigma}\left|\nabla(|x|^{-1} u)\right|^2~dx-(N-3)
\int_{\mathcal C_\Sigma}|x|^{-4}|u|^2~dx\\
&\ge&\left[\left(\frac{N-2}{2}\right)^2+\lambda_\Sigma-(N-3)\right]
\int_{\mathcal C_\Sigma}|x|^{-4}|u|^2~dx\\
&=&
\left[\left(\frac{N-4}{2}\right)^2+\lambda_\Sigma\right]
\int_{\mathcal C_\Sigma}|x|^{-4}|u|^2~dx~\!,
\end{eqnarray*}
as desired.
\QED

From now on we assume that 
\begin{equation}
\label{eq:cone}
\gamma_N+\lambda_\Sigma >0,
\end{equation}
where $\gamma_N=N(N-4)/4$, as before. Assumption (\ref{eq:cone})
is satisfied if $\Sigma=\Sf^{N-1}$ and $N\ge 5$.
In addition, the cases $N\ge 3$ and $\mathcal C_\Sigma=\R^N_+=$ homogeneous
half-space are included, as the first Dirichlet eigenvalue of the
Laplace-Beltrami operator on the half-sphere in $\R^N$ is $N-1$. In case $N=2$,
assumption (\ref{eq:cone}) is satisfied by any strictly convex cone.

If (\ref{eq:cone})
holds, then in particular the Hardy constant in inequality (\ref{eq:Hardy})
is positive, and
we can define a Hilbert space
$\mathcal D^{1,2}(\mathcal C_\Sigma;|x|^{-2}~dx)$ of maps $u$
vanishing on $\partial\mathcal C_\Sigma$ (if not empty), and such that
$$
\|u\|^2=\int_{\mathcal C_\Sigma}|x|^{-2}|\nabla u|^2~dx<\infty.
$$
Then $\mathcal D^{1,2}(\mathcal C_\Sigma;|x|^{-2}~dx)$ is continuously
embedded into
$L^2(\mathcal C_\Sigma;|x|^{-4}dx)$. Next we put
$$
\mathcal{N}^{2}(\mathcal{C}_{\Sigma})=
\left\{ u\in\mathcal D^{1,2}(\mathcal C_\Sigma;|x|^{-2}~dx)~|~
\Delta u\in L^2(\mathcal C_\Sigma)\right\}.
$$
By the Rellich inequality on cones proved in \cite{CM1} and using (\ref{eq:cone}) we get that
$$
\int_{\mathcal C_\Sigma}|\Delta u|^2~dx\ge
\left(\gamma_N+\lambda_\Sigma\right)^{2}~\!
\int_{\mathcal C_\Sigma}|x|^{-4}|u|^2~dx\quad\textrm{
for any $u\in \mathcal{N}^{2}(\mathcal{C}_{\Sigma})$.}
$$
Thus $\mathcal{N}^{2}(\mathcal{C}_{\Sigma})$ is a Hilbert space with norm
$$
\|u\|^2=\int_{\mathcal C_\Sigma}|\Delta u|^{2}dx,
$$
and $\mathcal{N}^{2}(\mathcal{C}_{\Sigma})\hookrightarrow
L^2(\mathcal C_\Sigma;|x|^{-4}dx)$.

A density argument can be used to show that $\mathcal{N}^{2}(\R^N\setminus\{0\})=\mathcal D^{2,2}(\R^N)$
if $N\ge 5$. Thus the next proposition includes Lemma \ref{L:regularity}
by taking $\Sigma=\Sf^{N-1}$.

\begin{Proposition}
\label{P:regularity}
Let $\Sigma$ be a domain of class $C^2$ in the unit sphere $\Sf^{N-1}$, and
assume that $\gamma_N+\lambda_\Sigma>0$. Then for any 
$g\in L^2(\mathcal C_\Sigma)$ there exists a unique
$v\in\mathcal{N}^{2}(\mathcal{C}_{\Sigma})$ such that
$$
-\Delta v=g\quad\textit{in $\mathcal C_\Sigma$.}
$$
Moreover, $v$ satisfies the Navier boundary conditions
$$
v=\Delta v=0\quad\textit{on $\partial\mathcal C_\Sigma$}
$$
if $\Sigma$ is properly contained in $\Sf^{N-1}$.
\end{Proposition}

\proof
We only have to prove existence. Uniqueness easily follows.
For any $R>1$ we put
$$
A^\Sigma_R=\left\{r\sigma~|~R^{-1}<r<R~,~\sigma\in\Sigma~\right\}.
$$
Let $v_R\in H^1_0(A^\Sigma_R)$ be the unique solution to
\begin{equation}\label{vRequation}
\begin{cases}
-\Delta v_R=g&\textrm{in $A^\Sigma_R$}\\
v=0&\textrm{on $\partial A^\Sigma_R$.}
\end{cases}
\end{equation}
We denote by $v_R\in \mathcal D^{1,2}(\mathcal C_\Sigma;|x|^{-4}dx)$
the null extension of $v_R$. Our aim is to 
use $|x|^{-2}v_R$ as test function in (\ref{vRequation}). Notice that
\begin{eqnarray*}
\int_{A^\Sigma_R}(-\Delta v_R)|x|^{-2}v_R&=&\int_{A^\Sigma_R}|x|^{-2}|\nabla v_R|^2~dx+
\frac{1}{2}\int_{A^\Sigma_R}\nabla (|x|^{-2})\cdot \nabla|v_R|^2~dx\\
&=&\int_{A^\Sigma_R}|x|^{-2}|\nabla v_R|^2~dx+(N-4)
\int_{A^\Sigma_R}|x|^{-4}|v_R|^2~dx
\end{eqnarray*}
and therefore from (\ref{vRequation}) we get
\begin{eqnarray}
\nonumber
&&
\left(\int_{A^\Sigma_R}|g|^2~\!dx\right)^{1/2}
\left(\int_{A^\Sigma_R}|x|^{-4}|v_R|^2~\!dx\right)^{1/2}
\\
\label{eq:CS}
&&\quad\quad \quad \quad \quad \quad\ge
\int_{A^\Sigma_R}|x|^{-2}|\nabla v_R|^2~dx+(N-4)
\int_{A^\Sigma_R}|x|^{-4}|v_R|^2~dx\\
\nonumber
&&\quad\quad \quad \quad \quad \quad\ge
(\gamma_N+\lambda_\Sigma)~\!\int_{A^\Sigma_R}|x|^{-4}| v_R|^2~dx
\end{eqnarray}
by the Cauchy-Schwarz inequality
and by Lemma \ref{L:Hardy}. Since $\gamma_N+\lambda_\Sigma>0$
by assumption, we first infer that $v_R$ is bounded in
$L^2(\mathcal C_\Sigma;|x|^{-4}dx)$. Thus, using (\ref{eq:CS})
again, we conclude that 
$v_R$ is uniformly bounded in $\mathcal D^{1,2}(\mathcal C_\Sigma;|x|^{-2}~dx)$.
Therefore,
up to a sequence $R\to\infty$, we can assume that $v_R\deb v$ weakly in
$\mathcal D^{1,2}(\mathcal C_\Sigma;|x|^{-2}~dx)$. It is easy to prove that
$v\in \mathcal D^{1,2}(\mathcal C_\Sigma;|x|^{-2}~dx)\hookrightarrow
L^2(\mathcal C_\Sigma;|x|^{-4}dx)$ solves $-\Delta v=g$ on $\mathcal C_\Sigma$. Thus in particular
$\Delta v=-g\in L^2(\mathcal C_\Sigma)$, that implies $v\in\mathcal N^2(\mathcal C_\Sigma)$. In case $\partial\mathcal C_\Sigma$ is not empty, then $v$ satisfies
Navier boundary conditions by standard arguments.
\QED

\medskip

\small
\noindent
{\bf Acknowledgements.} This research was done when the first author was visiting the
ICTP Mathematics Section in Trieste, Italy, and travel grant was sponsored by National Board of Higher Mathematics (NBHM), India. Warm hospitality of ICTP  is gratefully acknowledged. 

The second Author wishes to thank Prof. Fabio Zanolin for many helpful discussions
about equation (\ref{eq:ODE}) and for having suggested the references \cite{AT},
\cite{CT}.

\normalsize

\label{References}

\end{document}